\newcommand{\CLASSINPUTtoptextmargin}{1.8cm}
\newcommand{\CLASSINPUTbottomtextmargin}{3.6cm}
\documentclass[journal]{IEEEtran}
\usepackage{cite}
\usepackage{amsmath,amssymb,amsfonts}
\usepackage{algorithmic}
\usepackage[dvipdfmx]{graphicx}
\usepackage[font=footnotesize]{subcaption}
\usepackage{textcomp}
\usepackage{xcolor}
\usepackage{siunitx}
\usepackage{color,soul}
\usepackage{bm}
\usepackage{soul}
\usepackage{float}
\def\BibTeX{{\rm B\kern-.05em{\sc i\kern-.025em b}\kern-.08em
    T\kern-.1667em\lower.7ex\hbox{E}\kern-.125emX}}

\graphicspath{{../fig/}}

\DeclareCaptionFont{mysize}{\fontsize{8}{9.6}\selectfont}
\begin{document}
\bstctlcite{IEEEexample:BSTcontrol}

\title{Generalized Cancellation of Capacitor Parasitic Inductance Using a Lattice Network and Its Application to Common-Mode Noise Reduction}
\author{Katsuya~Nomura,~\IEEEmembership{Member,~IEEE,}~
        Shuhei~Chizuwa,~\IEEEmembership{}~
        Takashi~Masuzawa,~\IEEEmembership{Member,~IEEE}
\thanks{K. Nomura is with School of Engineering, Kwansei Gakuin University, Sanda 669-1330, Japan (e-mail:katsuya.nomura@kwansei.ac.jp).}
\thanks{S. Chizuwa is with Mitsubishi Heavy Industries, Yokohama 231-8715, Japan.}
\thanks{T. Masuzawa is with Mitsubishi Heavy Industries, Nagoya 453-8515, Japan.}
}

\maketitle

\begin{abstract}
This letter presents a generalized technique for cancelling the parasitic inductance of capacitors using a lattice network. 
{The $Z$-matrix-based equivalent transformation shows that the conventional vertically symmetric inductance condition can be relaxed for cancellation of the dominant high-frequency inductive component.
A vertically asymmetric lattice network is then proposed, enabling this component to be cancelled by adjusting only one inductance.}
The method is further extended to common-mode noise reduction using a common-mode choke, and both concepts are experimentally verified.
\end{abstract}


\section{Introduction}

Capacitors are essential components in electromagnetic compatibility (EMC) design. They are used to bypass conducted noise by utilizing their low-impedance characteristics at high frequencies. However, it is well known that the performance of practical capacitors is significantly degraded at high frequencies by equivalent series inductance (ESL), which represents the parasitic inductance of the capacitor~\cite{ott2011electromagnetic}.
Therefore, when bypass capacitors are mounted, it is important to reduce ESL, for example, by shortening the connection path~\cite{hubing1995power}. Nevertheless, ESL can be reduced but cannot be completely eliminated. In power converters, film capacitors and electrolytic capacitors are often used to satisfy requirements for large capacitance and high ripple-current capability. These capacitors tend to have larger ESL than small ceramic capacitors~\cite{wang2014reliability}.

As a technique that can theoretically eliminate the influence of ESL, ESL cancellation has been proposed in the field of power electronics. 
ESL cancellation techniques can be broadly classified into a method using mutual inductance between two inductors~\cite{Neugebauer2004,neugebauer2004filters} and a method forming a lattice network using two capacitors~\cite{wang2006cancellation}.
In the former method, among multiple inductive elements in a filter, the two inductors used for magnetic coupling can be selected with a certain degree of freedom~\cite{wang2024utilizing}, and several subsequent studies have been reported~\cite{mcdowell2013parasitic,kobayashi2022new}.
In contrast, the latter method realizes cancellation by the equivalent transformation of a lattice network rather than by magnetic coupling. However, subsequent studies on this approach have been relatively limited~\cite{wang2008study}. One possible reason is that the lattice-network-based method has been considered to require all ESLs in the circuit to be equal for ESL cancellation, resulting in a strong constraint and low implementation flexibility.

{In this letter, the condition for cancelling the dominant
high-frequency inductive component using a lattice network is
reconsidered, and it is shown that the conventional vertically
symmetric inductance condition can be relaxed. A Z-matrix-based
equivalent transformation reveals a more general condition under
the high-frequency approximation. Based on this result, an ESL
cancellation circuit using a vertically asymmetric lattice network
is proposed.}
The proposed method is further extended to common-mode
(CM) noise reduction using a CM choke.


The remainder of this letter is organized as follows. Section~II derives the ESL cancellation condition on the basis of the T-equivalent representation of a lattice network. Section~III experimentally verifies the effectiveness of the proposed method through $S$-parameter measurements using prototype circuits. Finally, Section~IV concludes this letter.

\section{Principle and Proposed Method}

\def\dd#1#2{\frac{\partial #1}{\partial #2}}
\def\muri{\mu_{r\_i}}

\begin{figure}[t]
\begin{minipage}[b]{0.5\linewidth}
\centering
\includegraphics[width=.8\linewidth,page=1]{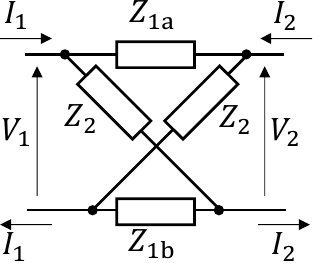}
\subcaption{}
\end{minipage}%
\hfil
\begin{minipage}[b]{0.5\linewidth}
\centering
\includegraphics[width=.8\linewidth,page=2]{fig-crop.pdf}
\subcaption{}
\end{minipage}%
\caption{(a) Lattice network and (b) T-equivalent circuit.}
\label{circuit}
\end{figure}

Fig.~\ref{circuit}(a) shows a lattice network. This circuit can be represented by the T-equivalent circuit shown in Fig.~\ref{circuit}(b), as demonstrated by comparing the $Z$-matrices of the two circuits.

First, the $Z$-matrix of the lattice network is derived. The diagonal
element is
\begin{align}
Z_{11}
&= \left. \frac{V_1}{I_1} \right|_{I_2=0}\notag \\
&= (Z_{1a}+Z_2) \parallel (Z_{1b}+Z_2)\notag \\
&=
\frac{(Z_{1a}+Z_2)(Z_{1b}+Z_2)}
     {Z_{1a}+Z_{1b}+2Z_2}.
\end{align}
Similarly,
\begin{align}
Z_{12}
&= \left. \frac{V_1}{I_2} \right|_{I_1=0} \notag \\
&=
\frac{(Z_{1b}+Z_2)Z_2 - (Z_{1a}+Z_2)Z_{1b}}
     {Z_{1a}+Z_{1b}+2Z_2} \notag \\
&=
\frac{Z_2^2 - Z_{1a}Z_{1b}}
     {Z_{1a}+Z_{1b}+2Z_2}.
\end{align}
{From the symmetry of the circuit, $Z_{21}=Z_{12}$ and
$Z_{22}=Z_{11}$.} Therefore, the $Z$-matrix $[Z_{\rm sym}]$ of the lattice network is given by
\begin{equation}
[Z_{\rm sym}] = 
\begin{bmatrix}
\dfrac{(Z_{1a}+Z_2)(Z_{1b}+Z_2)}{Z_{1a}+Z_{1b}+2Z_2} & \dfrac{Z_2^2 - Z_{1a}Z_{1b}}{Z_{1a}+Z_{1b}+2Z_2} \\
\dfrac{Z_2^2 - Z_{1a}Z_{1b}}{Z_{1a}+Z_{1b}+2Z_2} & \dfrac{(Z_{1a}+Z_2)(Z_{1b}+Z_2)}{Z_{1a}+Z_{1b}+2Z_2}
\end{bmatrix}.
\label{eq:sym}
\end{equation}

{The $Z$-matrix $[Z_{\rm T}]$ of the T-equivalent circuit is}
\begin{equation}
[Z_{\rm T}] = 
\begin{bmatrix}
Z_1' + Z_2' & Z_2' \\
Z_2' & Z_1' + Z_2'
\end{bmatrix}.
\label{eq:t}
\end{equation}

By comparing \eqref{eq:sym} and \eqref{eq:t}, the following relations are obtained:
\begin{align}
Z_1'
&=
\frac{
2Z_{1a}Z_{1b}
+
(Z_{1a}+Z_{1b})Z_2
}
{Z_{1a}+Z_{1b}+2Z_2}, \label{z1p}\\
Z_2'
&=
\frac{Z_2^2 - Z_{1a}Z_{1b}}
     {Z_{1a}+Z_{1b}+2Z_2}.
\label{z2p}
\end{align}
In particular, when $Z_1=Z_{1a}=Z_{1b}$ holds,
\begin{align}
Z_1'
&=
Z_1,\\
Z_2'
&=
\frac{Z_2 - Z_1}
     {2}
\label{z2p2}
\end{align}
is obtained.

\begin{figure}[t]
\centering
\includegraphics[width=.4\linewidth,page=3]{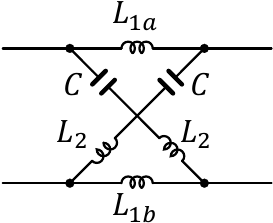}
\caption{ESL cancellation circuit.}
\label{cir_cancel}
\vspace{1mm}
\centering
\includegraphics[width=.4\linewidth,page=4]{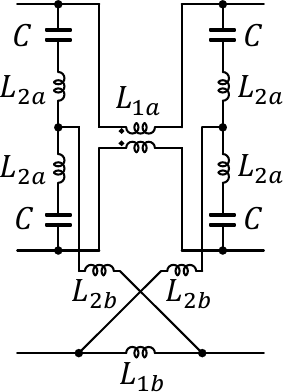}
\caption{ESL cancellation circuit for CM noise reduction.}
\label{cir_cancel_cm}
\end{figure}

Fig.~\ref{cir_cancel} shows the circuit diagram of ESL cancellation using a lattice network.
{As in the conventional method~\cite{wang2006cancellation},
mutual inductance and parasitic coupling between the branches are not
included in the present analytical model. Their quantitative effects
will be investigated in future work.}
In the conventional method~\cite{wang2006cancellation}, it was explained that ESL cancellation requires
\begin{align}
     {L_2}=L_{1a}=L_{1b}. \label{cancel_org}
\end{align}
Under this condition, the effective ESL of the capacitors is cancelled and becomes zero, as is clear from \eqref{z2p2}. However, this condition requires all four ESLs to be equal and therefore imposes a strong implementation constraint.
{At frequencies where the ESL impedance is dominant over the
capacitive reactance in the shunt paths, the inductive component of
$Z_2'$ in \eqref{z2p} can be approximated as
\begin{align}
Z'_{2,{\rm ind}}
&\simeq j\omega
\frac{L_2^2-L_{1a}L_{1b}}
     {L_{1a}+L_{1b}+2L_2}.
\end{align}
Under this approximation, the dominant high-frequency inductive component
can be cancelled even when $L_{1a}$ and $L_{1b}$ are unequal, provided that}
{
\begin{align}
     L_2^2=L_{1a}L_{1b} \label{cancel_prop}
\end{align}
is satisfied. Therefore, the conventional condition in
\eqref{cancel_org} is one sufficient condition, and the dominant
high-frequency inductive component can also be cancelled with a
vertically asymmetric lattice network.}

{If a single vertically asymmetric lattice network shown in
Fig.~\ref{cir_cancel} is directly connected between the positive (P) and
negative (N) conductors, the symmetry between the two conductors is broken,
which may cause conversion between CM and differential-mode (DM) signals.
To maintain the symmetry of
the overall two-line system, identical asymmetric lattice networks can instead
be connected separately between P and ground (GND) (the P--GND path) and
between N and GND (the N--GND path). In this configuration, the corresponding
element values and layouts should be matched between the two paths.}

{A single circuit configuration for CM noise reduction can also be
constructed while preserving the symmetry between the P and N conductors,
as shown in Fig.~\ref{cir_cancel_cm}.} In this configuration, a CM choke provides inductance that does not
affect the differential-mode. This inductance is not intended to directly
suppress the CM current through the high impedance of the choke, as in a
conventional EMI filter. Instead, it functions as an adjustment element for
satisfying the ESL cancellation condition of the lattice network.
{This CM noise reduction configuration cannot be realized in the
same manner under the equality constraint of the conventional method and is
enabled by the relaxed cancellation condition.}

\section{Experimental Verification}
\begin{figure}[tb]
\begin{minipage}[b]{0.99\linewidth}
\centering
\includegraphics[width=.5\linewidth,page=5]{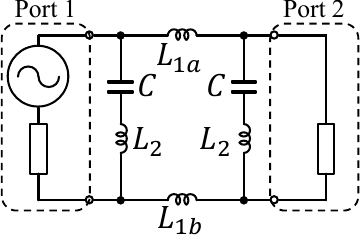}
\subcaption{}
\end{minipage}%
\hfil
\vspace{1mm}
\begin{minipage}[b]{0.99\linewidth}
\centering
\includegraphics[width=.7\linewidth,page=7]{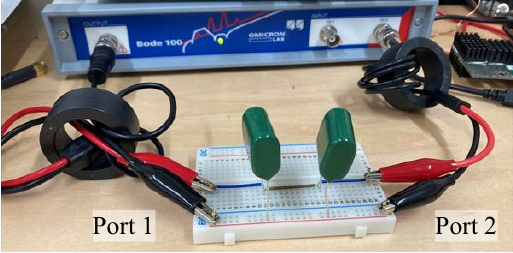}
\subcaption{}
\end{minipage}%
\caption{Normal circuit and measurement setup: (a) circuit diagram and (b) photograph.}
\label{board_normal}
\end{figure}

Prototype circuits were fabricated to verify the effectiveness
of the proposed method. Figs.~\ref{board_normal}(a) and (b)
show the normal circuit used for the measurement and a photograph
of the setup, respectively. The circuits were constructed using a
breadboard, jumper wires, and two capacitors (Rubycon,
450MPS105J, \SI{1}{\micro\farad}). The two-port S-parameter
$S_{21}$ was measured using a network analyzer (OMICRON Lab,
Bode 100). To suppress the influence of ground-loop-induced
CM currents, CM chokes were formed by
winding each port cable around nanocrystalline soft magnetic cores
(Proterial, FT-3K50T F6045GS) during the measurement~\cite{Davis2019GroundLoop}.

\begin{figure}[tb]
\centering
\includegraphics[width=.3\linewidth,page=6]{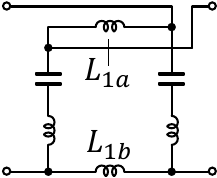}
\caption{Circuit diagram of the ESL cancellation circuit.}
\label{board_cancel_sch}
\vspace{1mm}
\begin{minipage}[b]{0.5\linewidth}
\centering
\includegraphics[width=.7\linewidth,page=8]{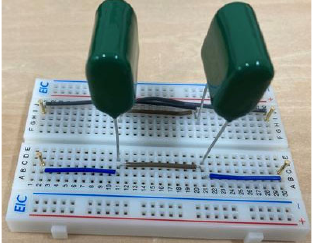}
\subcaption{}
\end{minipage}%
\hfil
\begin{minipage}[b]{0.5\linewidth}
\centering
\includegraphics[width=.67\linewidth,page=9]{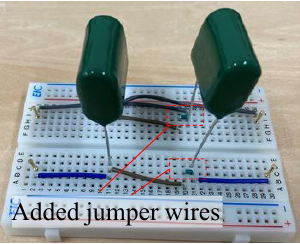}
\subcaption{}
\end{minipage}%
\caption{Conventional cancellation circuit: (a) without and (b) with additional jumper wires.}
\label{board_conv}
\vspace{1mm}
\begin{minipage}[b]{0.49\linewidth}
\centering
\includegraphics[width=.7\linewidth,page=10]{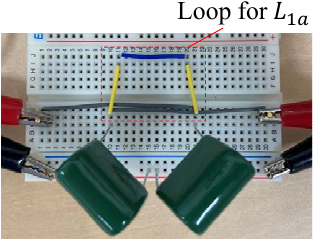}
\subcaption{}
\end{minipage}%
\hfil
\begin{minipage}[b]{0.33\linewidth}
\centering
\includegraphics[width=.99\linewidth,page=11]{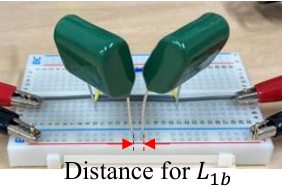}
\subcaption{}
\end{minipage}
\vspace{1mm}
\begin{minipage}[b]{0.32\linewidth}
\centering
\includegraphics[width=.99\linewidth,page=12]{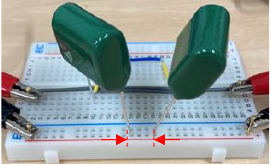}
\subcaption{}
\end{minipage}%
\hfil
\begin{minipage}[b]{0.32\linewidth}
\centering
\includegraphics[width=.99\linewidth,page=13]{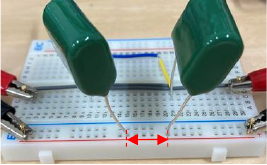}
\subcaption{}
\end{minipage}%
\hfil
\begin{minipage}[b]{0.33\linewidth}
\centering
\includegraphics[width=.99\linewidth,page=22]{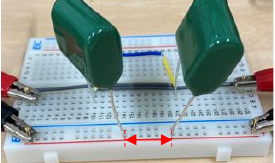}
\subcaption{}
\end{minipage}
\caption{Proposed cancellation circuit: {(a) loop
corresponding to $L_{1a}$; (b)--(e) loops with the conductor distance
corresponding to $L_{1b}$ set to \SI{2.54}{\milli\meter},
\SI{7.62}{\milli\meter}, \SI{12.70}{\milli\meter}, and
\SI{15.24}{\milli\meter}, respectively}.}
\label{board_prop}
\vspace{1mm}
\centering
\includegraphics[width=.82\linewidth,page=14]{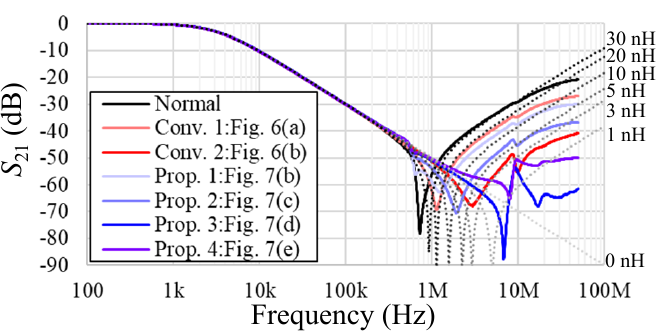}
\caption{Measured results of $S_{21}$.}
\label{s21}
\end{figure}

Fig.~\ref{board_cancel_sch} shows the basic ESL cancellation circuit.
Figs.~\ref{board_conv} and \ref{board_prop} show the circuits based on
the conventional and proposed conditions, respectively. In the conventional
circuit, $L_{1a}$ and $L_{1b}$ must be adjusted simultaneously to satisfy
\eqref{cancel_org}. 
In contrast, in the proposed circuit, a loop corresponding
to $L_{1a}$ is introduced, making $L_{1a}$ clearly larger than
$L_{1b}$. Then, $L_{1b}$ is adjusted while $L_{1a}$ is fixed by
changing the distance corresponding to $L_{1b}$. 
Thus, the proposed circuit can satisfy the cancellation condition in
\eqref{cancel_prop} by adjusting only one inductance.
{The generalized condition instead provides additional design
freedom when the two inductive paths are subject to different geometrical
or implementation constraints.}

Fig.~\ref{s21} shows the measured $S_{21}$ results. The dashed curves
represent simulated results of Fig.~\ref{board_normal}(a) for different
values of $L_2$ with $C=\SI{1}{\micro\farad}$ and
$L_{1a}=L_{1b}=\SI{0}{\nano\henry}$. The normal circuit exhibits an ESL
of approximately \SI{30}{\nano\henry}. Both the conventional and proposed
cancellation circuits reduce the effective ESL. 
In the proposed circuit, the best result is obtained for the third proposed
configuration shown in Fig.~\ref{board_prop}(d), where {the apparent
effective ESL estimated from the breadboard measurement is approximately
\SI{1}{\nano\henry} or less.}
The degradation observed when the conductor
spacing is either increased or decreased is consistent with \eqref{cancel_prop},
because $L_{1b}$ becomes either smaller or larger than the value required for
cancellation.
{The resonance and anti-resonance observed around
\SI{10}{\mega\hertz} may be associated with residual reactive components
outside the high-frequency approximation as well as implementation
parasitics. Their quantitative separation will be investigated in future
work.}

\begingroup
PCB prototypes with controlled conductor geometries were also evaluated to
reduce uncertainty due to the breadboard interconnections. Unlike the
breadboard prototypes, the capacitor leads were cut short before soldering.
Fig.~\ref{pcb_cancel} compares the normal, conventional, and proposed PCB
configurations. For each cancellation configuration, five PCBs were fabricated
with the conductor spacing varied in \SI{2}{\milli\meter} increments: from
\SI{16}{\milli\meter} to \SI{24}{\milli\meter} for the conventional
configuration and from \SI{6}{\milli\meter} to \SI{14}{\milli\meter} for the
proposed configuration.

\begin{figure}[!t]
\begin{minipage}[b]{0.48\linewidth}
\centering
\includegraphics[width=.99\linewidth,page=24]{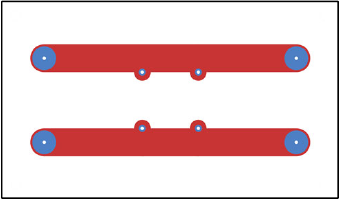}
\subcaption{Normal layout.}
\end{minipage}%
\hfil
\begin{minipage}[b]{0.48\linewidth}
\centering
\includegraphics[width=.99\linewidth,page=29]{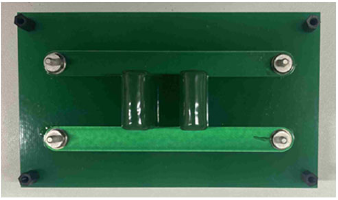}
\subcaption{Normal PCB.}
\end{minipage}
\vspace{1mm}
\begin{minipage}[b]{0.48\linewidth}
\centering
\includegraphics[width=.99\linewidth,page=27]{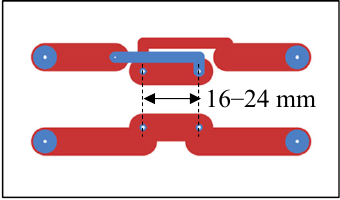}
\subcaption{Conventional layout.}
\end{minipage}%
\hfil
\begin{minipage}[b]{0.48\linewidth}
\centering
\includegraphics[width=.99\linewidth,page=30]{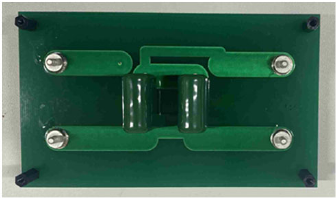}
\subcaption{Conventional PCB.}
\end{minipage}
\vspace{1mm}
\begin{minipage}[b]{0.48\linewidth}
\centering
\includegraphics[width=.99\linewidth,page=28]{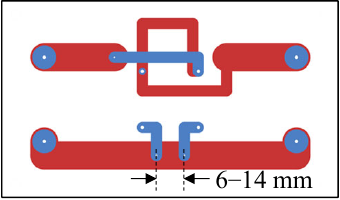}
\subcaption{Proposed layout.}
\end{minipage}%
\hfil
\begin{minipage}[b]{0.48\linewidth}
\centering
\includegraphics[width=.99\linewidth,page=31]{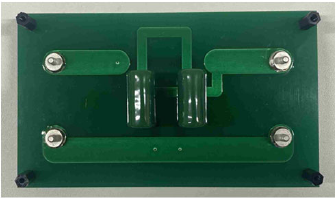}
\subcaption{Proposed PCB.}
\end{minipage}
\caption{PCB configurations used for ESL cancellation verification.}
\label{pcb_cancel}
\end{figure}

The measured results are shown in Fig.~\ref{s21_pcb}. The conventional and
proposed configurations both exhibit a nonmonotonic change as the conductor
spacing is varied, as in the breadboard results. The representative results
for the conventional spacing of \SI{22}{\milli\meter} and the proposed spacing
of \SI{10}{\milli\meter} confirm that both configurations reduce the influence
of ESL compared with the normal PCB.

\begin{figure}[!t]
\begin{minipage}[b]{0.99\linewidth}
\centering
\includegraphics[width=.9\linewidth,page=35]{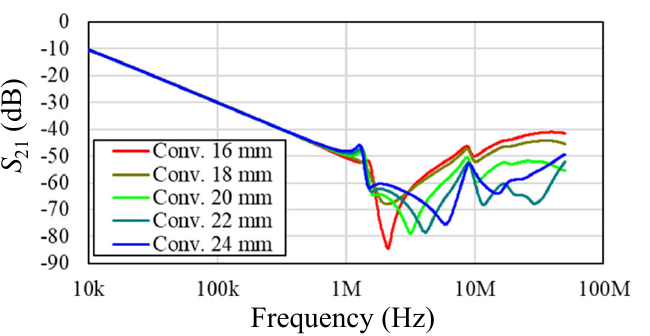}
\subcaption{Conventional configuration.}
\end{minipage}%
\vspace{1mm}
\begin{minipage}[b]{0.99\linewidth}
\centering
\includegraphics[width=.9\linewidth,page=36]{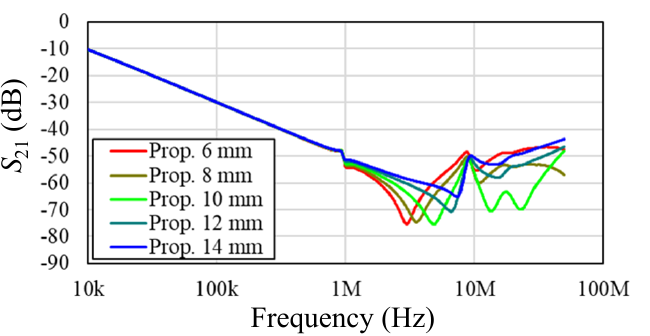}
\subcaption{Proposed configuration.}
\end{minipage}%
\vspace{1mm}
\begin{minipage}[b]{0.99\linewidth}
\centering
\includegraphics[width=.9\linewidth,page=37]{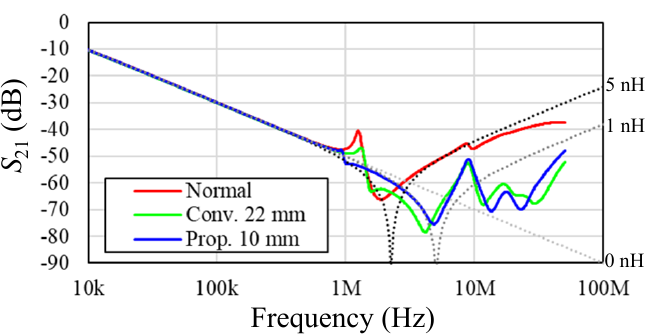}
\subcaption{Comparison of selected configurations.}
\end{minipage}
\caption{Measured $S_{21}$ results for the PCB verification.}
\label{s21_pcb}
\end{figure}
\endgroup

\begin{figure}[!t]
\begin{minipage}[b]{0.4\linewidth}
\centering
\includegraphics[width=.99\linewidth,page=20]{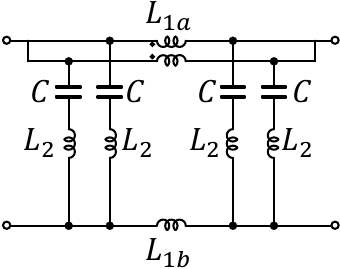}
\subcaption{}
\end{minipage}%
\hfil
\begin{minipage}[b]{0.4\linewidth}
\centering
\includegraphics[width=.99\linewidth,page=21]{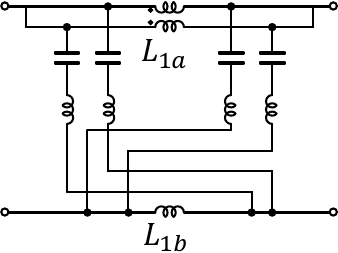}
\subcaption{}
\end{minipage}
\vspace{1mm}
\begin{minipage}[b]{0.52\linewidth}
\centering
\includegraphics[width=.99\linewidth,page=15]{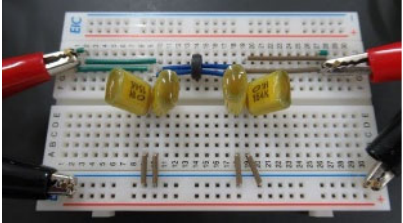}
\subcaption{}
\end{minipage}%
\hfil
\begin{minipage}[b]{0.45\linewidth}
\centering
\includegraphics[width=.99\linewidth,page=16]{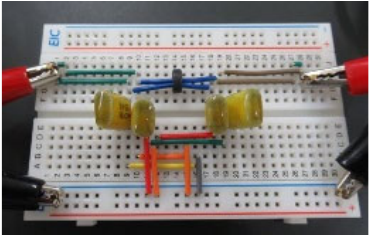}
\subcaption{}
\end{minipage}
\begin{minipage}[b]{\linewidth}
\centering
\includegraphics[width=.9\linewidth,page=19]{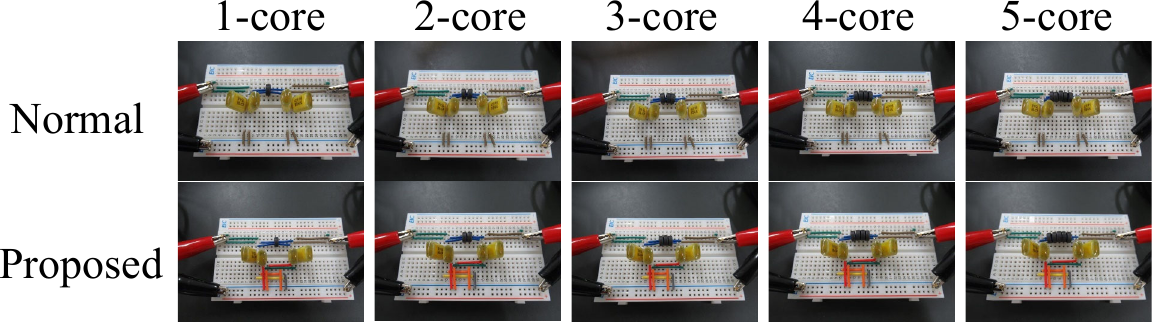}
\subcaption{}
\end{minipage}%
\caption{CM verification circuits: (a),(b) circuit diagrams and (c),(d) photographs of the normal and proposed circuits, respectively; (e) photographs of the configurations with different numbers of magnetic cores {(one to five cores from left to right, with the normal and proposed circuits in the upper and lower rows, respectively)}.}
\label{board_cm}
\end{figure}

\begin{figure}[!t]
\begin{minipage}[b]{0.99\linewidth}
\centering
\includegraphics[width=.82\linewidth,page=17]{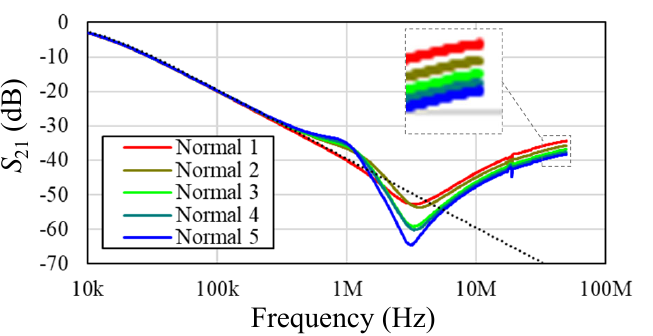}
\subcaption{}
\end{minipage}%
\hfil
\begin{minipage}[b]{0.99\linewidth}
\centering
\includegraphics[width=.82\linewidth,page=18]{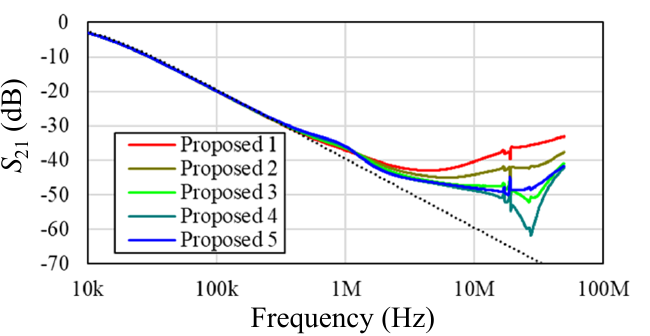}
\subcaption{}
\end{minipage}%
\caption{$S_{21}$ of the CM verification circuits: (a) normal circuit {including an enlarged view of the high-frequency region} and (b) proposed circuit.}
\label{s21_cm}
\end{figure}

The proposed method was also evaluated for CM noise reduction. Fig.~\ref{board_cm}
shows the normal and proposed CM verification circuits. Film capacitors
(Nissei Electric, AMZB0050J15400000000, \SI{150}{\nano\farad}) and Ni-Zn
ferrite cores (TDK, B64290P0037X001) were used, and $S_{21}$ was measured
while varying the number of cores corresponding to $L_{1a}$ from one to five.
{Because the cancellation effect was required to extend to several
tens of megahertz, Ni-Zn ferrite cores were selected for their high-frequency
characteristics instead of Mn-Zn ferrite cores.}
Fig.~\ref{s21_cm} shows the measured results. In the normal circuit, the $S_{21}$
changes almost monotonically with the number of cores. In contrast, 
the proposed circuit exhibits the lowest $S_{21}$ around the resonance region with four cores, and the $S_{21}$ degrades when the number of cores is either increased or decreased. This behavior agrees
with the tendency observed in Fig.~\ref{s21} and confirms that the proposed
cancellation concept is also effective for CM noise reduction.

\begingroup
The CM cancellation effect was also verified using PCBs. Fig.~\ref{pcb_cm}
shows the layouts and photographs of the normal and proposed configurations.
The crossing conductors in the proposed configuration are routed on different
copper layers and connected through vias.
During the PCB measurements, the P and N conductors were short-circuited using
copper plates, as shown in Fig.~\ref{pcb_cm}(d).

\begin{figure}[!t]
\begin{minipage}[b]{0.48\linewidth}
\centering
\includegraphics[width=.99\linewidth,page=38]{fig-crop.pdf}
\subcaption{Normal layout.}
\end{minipage}%
\hfil
\begin{minipage}[b]{0.48\linewidth}
\centering
\includegraphics[width=.99\linewidth,page=40]{fig-crop.pdf}
\subcaption{Normal PCB.}
\end{minipage}
\vspace{1mm}
\begin{minipage}[b]{0.48\linewidth}
\centering
\includegraphics[width=.99\linewidth,page=39]{fig-crop.pdf}
\subcaption{Proposed layout.}
\end{minipage}%
\hfil
\begin{minipage}[b]{0.48\linewidth}
\centering
\includegraphics[width=.99\linewidth]{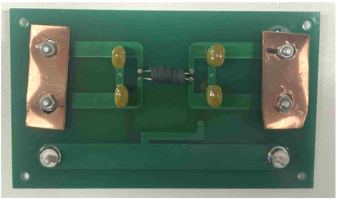}
\subcaption{Proposed PCB.}
\end{minipage}
\caption{PCB configurations for CM noise reduction.}
\label{pcb_cm}
\end{figure}

Fig.~\ref{s21_cm_pcb} shows the measured results obtained while varying the
number of magnetic cores from one to five. The proposed configuration exhibits
a distinct high-frequency minimum with three cores and changes when the number
of cores is either increased or decreased. Thus, the PCB results reproduce the
characteristic dependence on the adjusted inductance observed in the
breadboard experiment.

\begin{figure}[!t]
\begin{minipage}[b]{0.99\linewidth}
\centering
\includegraphics[width=.9\linewidth,page=43]{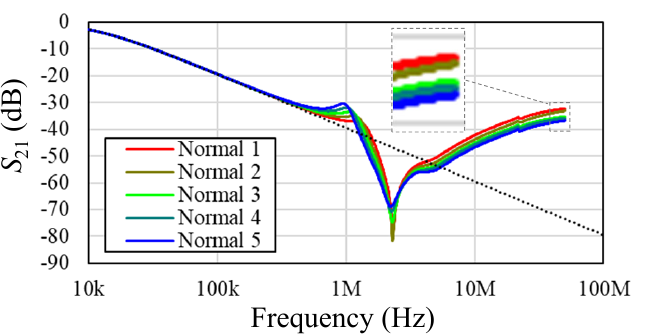}
\subcaption{Normal configuration.}
\end{minipage}%
\vspace{1mm}
\begin{minipage}[b]{0.99\linewidth}
\centering
\includegraphics[width=.9\linewidth,page=44]{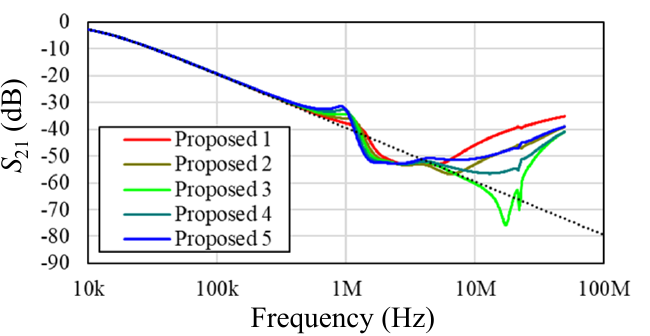}
\subcaption{Proposed configuration.}
\end{minipage}
\caption{Measured $S_{21}$ results for the CM PCB verification: (a) normal
configuration including an enlarged view of the high-frequency region and
(b) proposed configuration.}
\label{s21_cm_pcb}
\end{figure}

The required inductance is relatively small in the present small-signal
implementation. For high-power applications, magnetic-core saturation due to
CM excitation and imperfect cancellation of the DM flux must be considered.
\endgroup

\section{Conclusion}
{This letter generalized the condition for cancelling the dominant
high-frequency inductive component of capacitors using a lattice network. The
Z-matrix-based equivalent transformation showed that the conventional
vertically symmetric condition can be relaxed under the high-frequency
approximation. Based on this result, a vertically asymmetric cancellation
circuit was proposed, in which the condition can be satisfied by adjusting
only one inductance.} The method was also extended to CM noise reduction using
a CM choke.
The effectiveness of both the asymmetric ESL cancellation circuit and its CM-noise-reduction extension was verified through $S$-parameter measurements using prototype circuits.

\section*{Acknowledgment}
{The authors would like to thank Mr. Hayato Mitsuda for conducting
the PCB measurements.}

\bibliographystyle{IEEEtran}
\bibliography{IEEEabrv,reference_all}

\end{document}